
\magnification=\magstep1
\input xypic
\documentstyle{amsppt}

\def\noi{\noindent}
\def\Bigskip{\bigskip\bigskip}
\def\XT{{\Cal X}}
\def\Xs{{\Cal X_s}}
\def\X0{{\Cal X_0}}
\def\O{\Cal O}
\def\g{{\frak g}}
\def\h{{\frak h}}

\def\D{\Delta}

\define\pd#1#2{\dfrac{\partial#1}{\partial#2}}
\def\qed{\vrule height6pt width4pt}

\topmatter
\title DEGENERATION OF CALABI-YAU MANIFOLD  \\
WITH WEIL-PETERSSON METRIC
\endtitle
\author
Yoshiko Hayakawa
\endauthor
\rightheadtext{Degeneration of Calabi-Yau}
\address Department of Mathematics,
401 Mathematical Sciences,
Oklahoma State University,
Stillwater, OK 74078
\endaddress

\abstract
Koiso was first to introduce the Weil-Petersson metric in higher
dimension.
Tian showed that a moduli space of Calabi-Yau $n$-manifolds comes
naturally with Weil-Petersson metric. In this paper
we focus on determining for which degenerations  the central fibre is
at finite distance with respect to
 the Weil-Petersson metric. First we  obtain a
simple condition on the limiting mixed Hodge structure which is a necessary
and sufficient condition for finite Weil-Petersson distance to  the central
fibre.
 This issue has been raised in the Physics literature but not
 extensively analyzed there.
 Then we  combine the result with the canonical mixed Hodge structure
of the central fibre and obtain a simple cohomological
 necessary and sufficient condition for
the central fibre to be at  finite distance. As a corollary, we prove
that a central fibre with simple nodes is at finite distance.
\endabstract
\toc
\specialhead {} Introduction \endspecialhead
\head 1.  Degeneration with Weil-Petersson Metric \endhead
\head 2.  Limiting Mixed Hodge structure \endhead
\head 3.  Central fibre and Clemens-Schmid exact sequence \endhead
\head 4.  Examples \endhead
\endtoc
\endtopmatter

\document
\specialhead {} Introduction. \endspecialhead
\bigskip
Koiso \cite{Ks} first extended classical Weil-Petersson metric to
higher dimension. Then Tian showed that
 a moduli space of Calabi-Yau $n-$manifolds comes naturally  with
 the Weil-Petersson metric defined by using a generator of $H^{n,0}$
,which first
appeared in \cite{Ti}. The problem of compactifying a global moduli space is
not yet solved for higher dimensions, except for some examples in dimension
$3$ worked out with regards to
Mirror Symmetry. Even so-called degeneration problems were not solved  for
$n \geq 3$.  In this paper we complete the moduli space of Calabi-Yau
$n-$folds   with respect to the Weil-Petersson metric; i.e., to fill the
central
fibre of a degeneration  which is at finite distance with respect to
 the Weil-Petersson metric.
This constitutes a partial compactification of a moduli space of Calabi-Yau
$n-$folds.

A degeneration of $n-$dimensional varieties is a proper flat holomorphic
map from a variety $\XT$ to a disk $\D$ of relative dimension
$n$ such that the generic fibre is a smooth variety.  In the
case of Calabi-Yau $n-$folds, a one-parameter family of quintic $3-$folds
in ${\bold CP}^4 $ is a typical example. Any degeneration can be reduced
to a semi-stable degeneration  by a base change,
thanks to Mumford's semi-stable reduction
theorem \cite{Mum}.  A semi-stable degeneration allows only a reduced
divisor with normal crossings for a central fibre,
 which is among all singularities the easiest
to treat.  A typical example arises from a full blow-up of a complete
resolution of simple singularities such as simple nodal points.
In fact by a glorious result of Hironaka, any singularity can be
resolved into such a reduced variety with only normal crossings

For lower dimensional cases corresponding to ${\bold K}_{\Xs}=0$ ,the
results are already known.  In the case of the degeneration of complex tori
the genus will strictly go down, i.e., the central fibre has to become
a ${\bold CP}^1 $.  For ${\bold K3}$ surfaces, there are
$3$ cases that were worked out recently by \cite{PP} and
\cite{Kul}.
 Both results were obtained by  strictly Algebro-Geometric methods.

On the other hand, the deformation theory of complex structures was
originally developed by Kodaira and Spencer \cite{KS}.
Griffiths \cite{Gr1} started the construction of the period map and
later variation of Hodge structures
to study moduli spaces of complex structures with fixed underlying
topological manifolds.  For complex varieties with ${\bold K}_{\Xs}=0$,
it is known that the deformation space of  complex structures is
unobstructed, by the theorem of  Bogomolov-Tian-Todorov \cite{Ti}, \cite
{Tod}.  As a result,
one can look at the variation of Hodge structures associated to a given
moduli space problem, in particular, a degeneration.

 From a Hodge theoretic point of view, the study of a degeneration corresponds
to the study of (limiting) mixed Hodge structures, which is in a sense
a `degenerated' Hodge structure of generic fibres.
It should be noted that a `degenerated'  Hodge structure may no longer
be a Hodge structure.
Geometrically, a classifying space of Hodge structures is a
subset of a product of Grassmannians. Since the automorphism group
of this classifying space is transitive, it is also a homogeneous space.
Hence, a  representation of
$\pi_1(\D)$, usually called  the monodromy group, acts on the cohomology of
generic fibres, and in the case of a semi-stable
degeneration the monodromy operator is in fact unipotent by the
Monodromy Theorem.

Independent of a variation of Hodge structures, any (possibly singular)
 algebraic variety
carries a canonical mixed Hodge structure, as was first proved
by Deligne  \cite{DeH2}, \cite{DeH3} in rather general language.
So far, one of the few tools that works regardless of the dimension of a
generic
fibre is an exact sequence of cohomologies (and homologies) with
underlying mixed Hodge structures, called the Clemens-Schmid exact
sequence. This can connect the limiting mixed Hodge structure
of a generic fibre with the canonical mixed Hodge structure of the central
fibre.

It is usually impossible to compute a `limit' of Hodge structures of
generic fibres. Instead, we look at the {\it nilpotent orbit} of
a limiting mixed Hodge structure. Schmid \cite{Sch} proved that a nilpotent
orbit behaves asymptotically as a real `limit' of Hodge structures.
Hence one can compute the distance to the central fibre using
this nilpotent orbit.

We  rely on the geometry
of the classifying space, which, as a subspace of a product  of
Grassmannians,
 carries a canonical metric.
 Schmid's  nilpotent orbit
theorem refers to the distance by this metric.
However, the Weil-Petersson metric associated to any Calabi-Yau
moduli space is bounded by this canonical metric.

Our main result is stated as the following Theorem;
\smallskip
\noi {\bf Theorem   }  {\it  A semistable degeneration of a Calabi-Yau
$n$ dimensional compact complex variety is at  finite distance with
respect to the Weil-Petersson metric if and only if $H^{n,0}$ of one
of the components of the central fibre has positive rank.}
\smallskip
In secition $1$,
 we  define the variation of Hodge structures
associated to a given degeneration of Calabi-Yau $n-$folds and then
 define the Weil-Petersson metric and  show it is essentially
a direct summand of the canonical metric. Hence, the distance with
respect to the Weil-Petersson metric is always bounded above by the
distance with respect to the canonical metric. As the first step to
prove the main theorem,
we show the monodromy condition
on the  polarization of the limiting mixed Hodge structure of
finite Weil-Petersson  distance.  In section $2$ we
 investigate  finite distance
limiting mixed Hodge structures  in more detail
 and obtain a monodromy condition . In section
$3$, we turn to the central fibre and look at the spectral sequence
that is used to describe the canonical mixed Hodge structure of
a central fibre.  Finally,  we combine the two constructions
by means of the Clemens-Schmid exact sequence to obtain a rather simple
cohomology
condition on the central fibre as a criterion for  finite distance.  In section
$4$, we  illustrate the result with several examples.  In particular,
we
generalize a result appearing in the physics literature (eg. \cite{CGH})
 which says that a central fibre with simple nodes is at finite distance
 for dimension  $3$. We show that this holds for all dimensions
greater
than $2$.
\Bigskip
\head 1. Degeneration with Weil-Petersson Metric \endhead
\bigskip

\noi For our purposes, a Calabi-Yau manifold is a  complex projective manifold
of dimension $n$ with trivial canonical bundle.
We will look at the degeneration space of Calabi-Yau $n-$folds over the unit
disk $\D$.

 Facts and proofs stated in the rest of this section can be found in the text
\cite{GrTA} Chapter 1 and Chapter 6.
\smallskip
\noi {\bf Definition 1.1 } {\it A   degeneration}
is a proper flat holomorphic map
$\pi\:\XT\longrightarrow\ \D\ $ of relative dimension $n$ such that
$\Xs = \pi^{-1}(s)$ is a smooth complex variety for $s \neq 0$ , and $\XT$
is a K\"ahler manifold.
A degeneration is  {\it semistable} if the central fibre $\X0$ is  a  divisor
with
(global) normal crossings; locally $\pi$ is defined by
$$ s=x_1 x_2 x_3 \ldots x_k$$
 i.e., $\X0$ is written as a sum of irreducible components where each
component  $X_i$ is  smooth and all components intersect transeversely
each other.
\smallskip
We  further assume  generic fibres to be projective and the total
space to be pseudo-projective.

Let $f^{\ast}\:\Cal X^{\ast} \longrightarrow\ \Delta^{\ast}\ $ be the
restriction to the punctured disk. Fix a
smooth fiber $X_s$. Since $f^{\ast}\ $is a $C^{\infty}\ $fibration,
$\pi_1(\Delta^{\ast})\ $acts on the cohomology
$H^m(\Cal X_s)$. The map $T\:H^m(\Cal X_s)\longrightarrow\ H^m(\Cal X_s)\
$induced by the canonical generator of
$\pi_1(\Delta^{\ast})\ $ is called  the {\it Picard-Lefschetz transformation}.

\smallskip
\noi {\bf Theorem (Monodromy Theorem)  }
{\it  \roster
\item $T\ $is quasi-unipotent, with index of unipotency at most $m\ $. In other
words, there is
some $k\ $such that $(T^k-I)^{m+1}=0$.

\item If $f\:\Cal X\longrightarrow\Delta\ $is a semi-stable degeneration, then
$T\ $is unipotent (i.e., we take $k=1$).
\endroster}
\smallskip
\noi {\it Proof }   \cite{LA}
\smallskip
\noindent We are going to look at the variation of Hodge structures associated
to a semi-stable degeneration of Calabi-Yau $n$-folds.
Hence, by the Monodromy Theorem, we may assume $T\ $is unipotent.
\smallskip
\noi {\bf Definition 1.2  }   A {\it Hodge structure} of weight $k$, denoted
$\{ H_{\bold Z}, H^{p,q} \} $ , is  given by a lattice $H_{\bold Z}$ of finite
rank together with a decomposition on its complexification $H=H_{\bold
Z}\otimes
{\bold C}$:
$$H=\underset{p+q=k}\to\bigoplus H^{p,q}$$
\noi such that
$$H^{p,q}=\overline{H^{q,p}}.$$

Since each fibre $\Xs$ is a complex projective manifold, by {\it the Hodge
Decomposition Theorem}, the complex de Rham cohomology of $\Xs$ in each
dimension can be written as
$$H^k_{\bold DR}(\Xs, \bold C)=\underset{p+q=k}\to\bigoplus H^{p,q}(\Xs)$$
giving a Hodge structure.

Furthermore,  since $\XT$ is a pseudo-projective variety,
the  K\"ahler 1-1 form $\omega$ is rational and, in particular,
its restriction to each fibre is rational.

We define
$$ L \: H^m(\Xs, {\bold C}) \longrightarrow H^{m+2}(\Xs, {\bold C}) \qquad\text
{ by cup with}\quad
\omega,$$
and we define the primitive cohomology of $\Xs$ by
$$ P^m(\Xs,{\bold C}) := {\bold Ker} L^{n-m+1} .$$

 This  may be defined  a priori up to a monodromy action; however,
since the 1-1 form comes from  the restriction
of the cohomology of the total space, these are invariant under
monodromy.

If we take the $H_{\bold Z}$ lattice of Hodge structures as the intersection
 of the primitive cohomology and the integral cohomology, we can define a
bilinear form
$$Q \: H_{\bold Z}\times H_{\bold Z} \rightarrow {\bold Z},$$
called the {\it polarization} on the Hodge structure as follows:
$$Q(\phi , \psi ) =(-1)^{k(k-1)/2}\int_{\Xs}\phi\wedge\psi\wedge\omega^{n-k}
,$$
\noi where $k$ is the weight, $n$ the dimension of $\Xs$, and $\omega$
the K\"ahler form.

Now we can define a polarized Hodge structure.
\smallskip
\noi {\bf Definition 1.3  } {\it A polarized Hodge structure} of weight $k$,
denoted $\{ H_{\bold Z}, H^{p,q} , Q \} $, is given by a Hodge structure of
weight $k$ together with a bilinear form
$$Q \: H_{\bold Z}\times H_{\bold Z} \rightarrow {\bold Z},$$
\noi which is symmetric for $k$ even and skew-symmetric for $k$ odd,
satisfying the two {\bf Hodge-Riemann} bilinear relations:
\roster
\item$  Q(H^{p,q}, H^{p',q'})=0 \quad\text {unless}\quad p'=k-p\quad
\text {and}\quad q'=k-q,$
\item$(\sqrt{-1} )^{p-q}Q(\psi,\overline\psi) > 0\quad\text {for
any nonzero element}\quad\psi\quad\text {in}\quad H^{p,q}.$
\endroster
\noi $Q$ is exactly  the polarization.

\smallskip
\noi {\it Remark 1.4  }  Equivalently, we can replace
the Hodge decomposition  by
a Hodge filtration, because a filtration varies holomorphicaly while
a decomposition does not . Given a decomposition $H=H^{p,q}$,
$$F^p=H^{k,0}\oplus\ldots\oplus H^{k-p,p}, \quad\text {and}\quad
H^{p,q}=F^P\cap\overline{F}^q.$$
\noi $\{ F^P \}$ is called the  {\bf Hodge filtration} of $H$.
With this notation, the polarization satisfies:
\roster
\item $Q(F^p,F^{k-p+1})=0$
\item $ Q(C\psi, \overline{\psi})>0 $ for any nonzero element
$ \psi$ in  $H$, where $C$ is the Weil operator of $H$ defined
by $C|_{H^{p,q}} = (\sqrt {-1})^{p-q}.$

\endroster

\smallskip
Hence, we have a polarized Hodge structure given by the Hodge decomposition
of each de Rham cohomology  group $H^k_{\bold DR}(\Xs, \bold C)$ with the
canonical
bilinear form as a polarization.  From now on, we
focus on the polarized Hodge structure of weight $n$, i.e.,
$H^n_{\bold DR}(\Xs, \bold C)$. Since each $\Xs$ is a Calabi-Yau $n$-fold,
the dimension of $H^{n,0}(\Xs)$ is always $1$. Furthermore the $H^{n,0}(\Xs)$
 part of $H^n_{\bold DR}(\Xs, \bold C)$ always carries the polarization
since $H^{n,0}(\Xs) \subset P^n(\Xs),$ which can  easily be seen
 to be primitive by the
construction of the usual Hodge diamond.

We would like to see how the polarized Hodge structure of each fibre varies.
Thus we need the notion of  a variation of (polarized) Hodge
structure. Here we use the construction done in  \cite{GrTA} Chapter 1.
Essentially a variation of Hodge structures is a
 holomorphic map from a moduli
space $S$ to a classifying space $D$, which is geometrically the open
set of the subset of the product of Grassmannians, and the map has
to satisfy certain tangential conditions.

The variation of Hodge structures
 of weight $n$ associated to the Calabi-Yau degeneration is defined
as the map $\phi$ in the following commutaive
diagram where $\h$ denotes the upper half plane that is the
universal cover of $\D^*$;
$$
\CD
\h   @>\tilde \phi>>    D \\
@V\pi VV                       @VVV \\
\D^*       @>\phi>>           D/ \Gamma
\endCD $$
\noi  and $\phi$ is defined by associating
to each $\Xs \ $ the Hodge decomposition of
$H^n_{DR}(\Xs,C)\ $. This is  the well-defined  map from $\D^*\ $
to  $D\ $up to monodromy action by $\Gamma$.
Since
$$\pd {H_s^{p,q}}{\bar s} \subseteq H^{p,q}\oplus H^{p+1,q-1}$$

\noi We have
$$\pd {F_s^p}{\bar s} \subseteq F_s^p$$
$$\text {and, by conjugation,}\quad\pd {F_s^p}s \subseteq F_s^{p-1}$$

We can have such a commutative diagram since we are assuming
that the degeneration is semi-stable and hence the monodromy is
unipotent.

\noi The polarization of weight $n$ variation of
Hodge structures becomes simplified as following;
$$Q(\phi , \psi ) =(-1)^{n(n-1)/2}\int_{\Xs}\phi\wedge\psi$$
\smallskip

\noindent {\it Remark 1.5 }   Let $M$ be a Calabi-Yau
 manifold. Then the map from the moduli space of $M$ to $D$
 defined by the  variation of Hodge structure is an immersion.
Since by deformation theory, the space of infinitesimal deformation of complex
structure
is given by $H^1(\Cal X_s,\Theta_s)\ $
and it is unobstructed by the theorem of Bogomolov-Tian-Todorov \cite{Ti},
\cite{Tod}.
\smallskip
Before we define the Weil-Petersson Metric
 we have to look at the  canonical metric on the classifying space
$D$.  We will call
this metric the VHS metric.

\noi For $\tilde\phi : \h \longrightarrow D\subset \Pi_{p=1}^nG(F^p,H)$,
we have  its differential
$$d\tilde\phi:T(\h ) \longrightarrow T(D)\subset T(\Pi_{p=1}^nG(F^p,H)).$$
Since this takes values in the horizontal subspace,
$$d\tilde\phi T(\h )\subset \underset {p=1}\to{\overset n\to\bigoplus} Hom(\Cal
H^{p,q},\Cal H^{p-1,q+1}).$$
In other words, we have a following commutative diagram,
$$
\CD
d\tilde\phi T(\h )    @>d\tilde\phi>>  \underset {p=1}\to{\overset
n\to\bigoplus} Hom(\Cal F^p / \Cal F^{p+1}, \Cal F^{p-1} / \Cal F^p ) \\
@ViVV                            @ViVV\\
T(D)       @>\subset>>  \underset {p=1}\to{\overset n\to\bigoplus} Hom(\Cal F^p
,\Cal H
 / \Cal F^p )
\endCD
.$$

\noi $H^{p,q}$ has a canonical Hermitian metric,  called
the Hodge metric, given by
$$<\psi,\eta>_{p,q}=(\sqrt {-1})^{p-q}Q(\psi,\bar \eta), $$
\noi where $Q(\ ,\overline\ )$ is the one we defined in section 1.

The canonical metric on $Hom(H^{p,q},H^{p-1,q+1}) \owns A$ is given by,

$$tr(A A^{\dagger})$$
where $A^{\dagger}$ is  the Hermitian conjugate such that
$$<A\alpha,\beta>_{p-1,q+1} =<\alpha,A^{\dagger}\beta>_{p,q}\qquad\text
{for}\quad
\alpha\in H^{p,q}, \beta\in H^{p-1,q+1}.$$

\noi Hence for
$$di:T(D)\longrightarrow \underset {p=1}\to{\overset n\to\bigoplus} Hom(\Cal
H^{p,q},\Cal H^{p-1,q+1}),$$
(where $ i$ is the inclusion  $ D\subset \Pi_{p=1}^nG(F^p,H)$,)
 and $v \in T(D)$, we have
$$
\split
<v,v>_{vhs}
 &=\underset {p=1}\to{\overset n\to\bigoplus}<d\phi(v),d\phi(v)>_{Hom_p}\\
 &=\sum_{p=1}^n trA_{v_p} A_{v_p}^{\dagger}\\ \endsplit
.$$
\noi Since dim$H^{n,0}=1$,  the summand for $p=n$ is just

$$tr(A_{v_n}A_{v_n}^{\dagger})= \frac {<A_{v_n}A_{v_n}^{\dagger}\Omega,\Omega>_
{n,0}}{<\Omega, \Omega>_{n,0} }= \frac {<A_{v_n}\Omega,A_{v_n}\Omega>_{n-1,1}}
{<\Omega, \Omega>_{n,0} }
.$$
where $0 \neq \Omega \in H^{n,0}$.

Now we define the Weil-Peterson Metric.  The classic Weil-Petersson
metric was only defined on a moduli space of curves.
Tian in \cite{Ti} extended the definition to moduli spaces of
Calabi-Yau $n-$folds.

Take a generator $w\in H_s^{n,0}(\Cal X_s,C)$. We know $H^{n,0}$ is a trivial
line bundle and it has a natural Hermitian metric
$(\sqrt{-1})^n Q(\,,\,)$ that satisfies the Hodge Riemman bilinear relation.

\noi We define
the Weil-Petersson metric to be a K\"ahler metric whose K\"ahler
form is $\frac 12 \Theta (s)$, i.e.,
$$ R(s)=-\bar\partial s\partial s log((\sqrt {-1})^n Q(w,\bar w))dsd\bar s.$$
The definition does not depend on a choice of
 $w\in H_s^{n,0}(\Cal X_s,C)$, since any  $w'\in H_s^{n,0}(\Cal X_s,C)$
is given by $ w'= hw $ where $h$ is a holomorphic function.

The proof  that it is indeed a metric  on the moduli
space is due to Tian
 \cite{Ti}.
\smallskip
\noi {\bf Proposition 1.6  } {\it The Weil-Peterson metric  is a  direct
summand of the VHS-metric. }
\smallskip
\noi {\it Proof  } See \cite{Ti}.

\smallskip
\noi {\bf Corollary 1.7  } {\it  Since Weil-Peterson metric is essentially a
summand of VHS metric, $\rho_{wp}$ Weil-Petersson distance is always bounded by
$\rho_D$ VHS metric
distance.}
\smallskip
\noi {\it Proof  }  It is clear by  inequality of integrands of each distance
integral.
\smallskip

Now we are ready to compute the distance  associated to a
degeneration of Calabi-Yau $n$-folds with respect to the Weil-Petersson
metric. In order to do this, we need the nilpotent orbit  originally
defined on  $\check D \supset D$, which
behaves asymptotically with the limit of Hodge filtrations
on $D$.
The construction of the nilpotent orbit is  given in
 Chapter $4$ of \cite{GrTA}.

Let $\phi \: \D^* \longrightarrow  D/ \Gamma $ be the variation of
Hodge structures associated to the family $\XT$ of Calabi-Yau three-folds
over the punctured disk.
Since $\Gamma$ has a representation in $\pi_1(\D^*)= {\bold Z}$, we may call
the
generator $T$. By  the Monodromy Theorem, $(T-I)^{n+1}=0$.
Let $\h=\{ w=u+vi | v>0 \}$ be the upper half plane. Then we
have a  commutative diagram:

$$
\CD
\h     @>\tilde\phi>>     D\\
@V\pi VV                        @VVV \\
\D^{*} @>\phi>> D/ \Gamma
\endCD
$$

\noi where,
$$ s=\pi(w)=\exp (2\pi iw)\quad\text
{and}\quad\tilde\phi(w+1)=T\tilde\phi(w).$$

\noi Analytic continuation around $s=0$ gives $H_{e^{2\pi is}}=TH_s$. Then

$$N:=\log T=(T-I)-\frac{(T-I)^2}{2} +\frac{(T-I)^3}{3} \ldots
\frac{(T-I)^n}{n}$$
Define
$$\tilde\psi:h\longrightarrow \check D\supset D \qquad \text {by} \qquad
\tilde\psi(w)=\exp (-wN)\tilde\phi(w)$$
Then
$$
\split
\tilde\psi(w+1)
 &=\exp (-(w+1)N)\tilde\phi(w+1) \\
 &=\exp (-wN)\cdot \exp (-N)\cdot T\tilde\phi(w) \\
 &=\exp (-wN)\cdot T^{-1}\cdot T\tilde\phi(w)  \\
 &=\exp (-wN)\tilde\phi(w)\\
 &=\tilde\psi(w)\endsplit
$$

So, $\tilde\psi$ descends to a single-valued map $\psi :\Delta^{\ast}
\longrightarrow \check D$ given by
$$
\split
\psi(s)
 &=\tilde\psi(w) \\
 &=\exp (\frac {-\log s}{2\pi i}N)\phi(s)\endsplit
.$$

Unlike $\phi$, the map $\psi$ extends over $s=0$ due to the following
Theorem.
\smallskip
\noi {\bf Theorem (Cornalba and Griffiths )} {\it The map $\psi:\Delta^{\ast}
\longrightarrow \check D$ extends across the origin to a map $\psi:\Delta
\longrightarrow \check D$ .}
\smallskip
\noi {\it Proof } See \cite{GrCo} or \cite{GrSc}.
\smallskip
\noi In fact, the $\psi (0)$ has a special name.
\smallskip
\noi {\bf Definition 1.8 }    The filtration $\psi(0)\in\check D$ will be
called the {\it limiting Hodge filtration} and will be denoted by
$F^P_{\infty}$.
\smallskip
\noi With the limiting Hodge filtration, we are ready to define
the nilpotent orbit.
\smallskip
\noi {\bf Definition 1.9  } {\it The nilpotent orbit} of a degenerating family
over $\Delta^{\ast}$ is the map $\Cal O:h \longrightarrow \check D$ given by
$\Cal O(w)=\exp (wN)\psi(0)$.
\smallskip
\noi Then
$$
\split
\Cal O(w+1) &=\exp ((w+1)N)\psi(0) \\
 &=\exp (wN)\exp (N)\psi(0) \\
 &=T\Cal O(w)\endsplit
$$
\noi Schmid showed in the following theorem that the nilpotent orbit
is a very good approximation of the period map.
\smallskip
\noi {\bf Theorem (Nilpotent orbit theorem )}
{\it
\roster
\item The nilpotent orbit is horizontal.
\item There is a non-negative number $\alpha $ such that, if $im(w) > \alpha $
, then $\Cal O(w)$ belong to $D$;
\item $\Cal O(w)$ is strongly asymptotic to $\tilde\phi(w)$ in the sense that
$\rho_D(\tilde\phi(w),\Cal O(w))\leq(Imw)^B e^{-2 \pi Imw}$
  for some $B\geq 0$ and $Imw\geq A > 0$ ,where $\rho_D$ is the distance on
  $D(\subset \check D)$ given be natural Hermetian Metric .
\endroster}
\smallskip
\noi {\it Proof } See \cite{Sch}.
\smallskip
\noi By Corollary 1.7 we immediately obtain
Weil-Petersson version of Nilpotent orbit Theorem. Furthermore,
even stronger version is true for Weil-Petersson metric.
\smallskip
\noi {\bf Lemma 1.10  } {\it The  central fibre of  a degeneration
has finite Weil-Petersson metric distance if and only if the Weil-Petersson
metric distance along the nilpotent orbit is finite.}
\smallskip

Before we prove the lemme, first we state and prove th following result.

\smallskip
\noi {\bf Theorem 1.11  } {\it Suppose the limit filtration $\{ F_{\infty}\}$
is given. Then the Weil-Peterson Metric distance of the corresponding nilpotent
orbit is
either $0$ or infinite.

Furthermore the Weil-Petersson Metric distance of the corresponding nilpotent
orbit is $0$, if and only if
$Q(\alpha,N^i \bar\alpha) =0 $ for all $i>0$ where $\alpha \in F_{\infty}^n$.}
\smallskip
\noi {\it Proof  } Let $\Cal O:h \longrightarrow \check D$ be the nilpotent
orbit, where
$$\Cal O(w)= e^{wN}F_{\infty}\quad\text {and}\quad \Cal O(w+1)=T\Cal O(w).$$
Since we want to compute the Weil-Petersson Metric distance of the nilpotent
orbit, we choose $\alpha \in F_{\infty}^n $ and set $\Omega(w)=e^{Nw} \alpha$.
Then,
$$R(w)=- \bar\partial w \partial w \log((\sqrt {-1})^n
Q(\Omega(w),\bar\Omega(w))$$
We look at the argument of the logarithm and get
$$
\split
(i)^n Q(\Omega(w),\bar\Omega(w))
 &=(i)^n Q(e^{Nw}\alpha,e^{N\bar w}\bar\alpha)\quad\text {since }\quad N
\quad\text {is real}\\
 &=(i)^n Q(e^{N(w_0+iy)}\alpha,e^{N(w_0-iy)}\bar\alpha)\\
 &=(i)^n Q(e^{Niy}\alpha,e^{-Niy}\bar\alpha)\ \text {since}\  e^{Nw_0}
 \ \text{fixes the  polarization}\\
 &=(i)^n Q(\alpha,e^{-2iyN}\bar\alpha)=(*)
\endsplit
$$
By the monodromy theorem, $N^{n+1}=0$. Hence
$$e^{-2iyN}=1-2iyN-2y^2N^2+\frac83iy^3N^3+ \ldots \frac{(-2i)^n}{n!}y^n N^n$$

\noi Then
$$
\split
(*)
 &=(i)^n Q(\alpha,(1-2iyN-2y^2N^2+\frac83iy^3N^3+\ldots +\frac{(-2i)^n}{n!}y^n
N^n)\bar\alpha)\\
 &=(i)^n Q(\alpha,\bar\alpha)+(i)^n Q(\alpha,-2iyN\bar\alpha)-\ldots +(i)^n
Q(\alpha, \frac{(-2i)^n}{n!} y^n N^n\bar\alpha)\endsplit
$$
Let $$(*)=p(y).$$
Then
$$R(w)=- (\frac {\partial^2}{\partial x^2} + \frac {\partial^2}{\partial
y^2})\log p(y).$$
Hence, the Weil-Petersson metric distance of the nilpotent orbit is given by
$$\rho_{wp}= \int_{y_0}^y \sqrt{- \frac {\partial^2}{\partial y^2}\log p(y)}
 dy$$
By  the Nilpotent Orbit Theorem, for $y >> 0$, $\Cal O  (w)\in D$.
 It follows that
$$(i)^n Q(\Omega(w),\bar\Omega(w)) > 0,$$
and hence $p(y)$ has to be positive.
Then the imaginary part of the polynomial vanishes for $y  >> 0$, i.e., it has
to be identically zero.

Therefore, $p(y)$ is  real and the leading coefficient is positive.
If we look at the leading term of $p(y) $, that is, $ay^k$, the argument of
the square root of the integrand becomes

$$\frac {-a^2(k-1)k y^{2k-2} + a^2 (k-1)^2y^{2k-2}}{a^2 y^{2k-2}}=
\frac {a^2}{y^2}$$
As $y$ approaches infinity, we may concentrate the effect of the
leading term. The distance is approximated by
$$\rho_{wp}= \lim_{y\to 0}\int_{y_0}^y \sqrt{\frac {a^2}{y^2}}
 dy$$
Hence the integral will be either infinite or zero, that is,  when
$a=0$.  The same argument works for the lower terms.
Hence, $\rho_{wp}$ is either infinite when all coefficients
are zero.
Now we go back to the definition of $p(y)$.
If we let
$$Q(\alpha,\bar\alpha)=C_0,\,Q(\alpha,N\bar\alpha)=C_1,\,Q(\alpha,N^2\bar\alpha)=C_2,\,\ldots \, Q(\alpha,N^n\bar\alpha)=C_n,$$
then this coefficient condition is true  if and only if
$$ C_1=C_2=C_3= \ldots =C_n=0.$$
 Hence the Weil-Petersson distance is finite if and only if
$$Q(\alpha,N^i\bar\alpha)=0 \forall i>0.$$
\line{\hfill \qed}
\smallskip
It remains to prove the {\bf Lemma 1.10}.
\smallskip
\noi {\it Proof of Lemma 1.10}
 Recall the map $\tilde\psi (w) =
\exp(-wN) (\tilde\phi) (w) $
 we described when we discussed  the construction of the nilpotent
orbit associated to the period map of the variation of Hodge
Structure:
$$
\CD
\h     @>\tilde\phi>>     D\\
@V\pi VV                        @VVV \\
\D^{*} @>\phi>> D/T
\endCD
$$
$\tilde\psi (w)$ is periodic with period $2\pi i$ and it descends to
the well-defined holomorphic map $\psi (s)$ from $\D^*$ to $\hat D$.
By the theorem of Cornalba and Griffiths, this map extends across the
origin and defined on $\D $.

In particular, $\Omega (s) = \exp (-wN)\hat\Omega = \exp(-\frac {\log
s}{2\pi i} N) \hat\Omega (s) $ is a single-valued holomorphic section
of the deepest level of the filteration $\exp (-wN) \tilde\phi (w)$
 that extends over $\D$, where $\hat\Omega $ is a constant global flat
section
and $\hat\Omega(s)$ is the push-foward that is multi-valued.
We may call $\Omega (s)$ a privileged section of the deepest level
of the filtration.

Choose  a frame $\{ \sigma_i \}$ of  privileged section  of
the Hodge bundle associated to the variation of Hodge structure.
Since it is given by the horizontal displacement of a basis of
the reference fibre, their values are constant.
$\Omega (s)$ can be written as a convergent power series around the
origin, of which coefficients $\Omega_k (s)$ are all privileded
sections:
$$
\split
\Omega (s)
 &= \sum_{k=0}^{\infty} \Omega_k  s^k \\
 &= \sum_{k=0}^{\infty} ( \sum_{i=1}^{n} w_{ki} \sigma_i ) s^k \qquad
\text{where} \quad w_{ki} \quad \text{are constant matrices}\\
 &= \sum_{i=1}^{n} ( \sum_{k=0}^{\infty} w_{ki} s^k ) \sigma_i \endsplit
$$
We choose a norm that is given by the following norm:
$$< \sigma_i, \sigma_j > = \delta_{ij} $$
Notice with this norm $ \| \Omega_k (s) \| $ is constant. Furthermore,
by  a  fact in basic linear algebra,
$Q(v,w) <  K \| v \| \| w \| =  K \sqrt {< v,v> <w,w> }$ for some
$K$. Let $\hat\sigma_i $'s be the multiple-valued descent of
the global flat section that is used to define a frame of
priviledged sections, i.e., $\sigma_i  = \exp(-wN) \hat\sigma_i $.
Then  since $Q(\  , \  )$ is monodromy invariant
$Q(\sigma_i ,  \sigma_j ) =Q(\hat\sigma_i , \hat\sigma_j) $
  ,where
$\hat\sigma_i$ are the pull-back  to the half-plane
$\h $ from $\D^* $ by the abuse of notation. Hence
$Q(\sigma_i, \sigma_j)$ is constant globally. Thus  the bound on
$Q(\ ,\  )$ by $\|  \  \| $ is uniform within the radious of
convergence.

By the abuse of notation, $\Omega (s)$ pull-back to,
$$
\split
\Omega (s)
 &= \sum_{k=0}^{\infty} \Omega_k  s^k \\
 &= \sum_{k=0}^{\infty} \Omega_k  \exp (2\pi i k w) \endsplit
$$
where $\Omega_k (w) $ is well-defind global section that is
the pull-back of the privileded section $\Omega_k $.

Then the deepest filtration of $\tilde\phi (w)$ , that is
$\Omega (w)$ can be written as:
$$ \Omega (w) = \exp (wN) \Omega (s) = \sum_{k=0}^{\infty}
\exp (wN) \Omega_k \exp (2\pi i k w). $$
Note that $\| \Omega_k \| $ is still constant, where
$\|\  \| $ is the pull-back norm of the one we defined.
Hence $\| \Omega_k \| \leq A^k $ for some $A>0$.

Finally since the nilpotent orbit map is given by
$\O (w) = \exp (wn)  \tilde\psi (0)$,
it's deepest filtaration part, which by the abuse of notation
denoted by $\O (w)$ is exactly,
$$\O (w)= \exp (wN) \Omega_0 $$
that is the constant term of $\Omega (w)$.
Let
$$\align
&\epsilon (w) = \Omega (w) - \O (w) \\
&\delta (w) = \frac {\partial \Omega }{\partial w} - \frac {\partial \O }
{\partial w}
\endalign$$
In order to estimate the differenc of two distance, first we estimate
the bounds of $\| \epsilon (w) \| $ and $\| \delta (w) \|$. We write
$$\epsilon (w) = \sum_{k=1}^{\infty} \exp ( 2 \pi i k w) \exp (Nw)
\Omega_k
  $$
Then along a vertical ray,
 $$\exp (Nw) (\Omega_k ) = \exp (N (w_0 + i y )) (\Omega_k )
 =\exp (Niy) (\Omega_k ) $$
$$\quad\text{since}\  Nw_0 \ \text{fixes} \  \Omega_k.$$
And $N^{n+1} = 0 $ gives;
$$\exp (Nw) (\Omega_k (w)) = ( 1 + Niy + \dots + \frac { (N i
y)^n}{n!} )
\Omega_k (w).$$
Finally, since there exist $B>0$ , which only depends on the choice of norm,
 such that $\| \Omega_k \| \leq B$ for all $k \geq 1 $ and $N $ raises
 norms by the bounded factor for all $k \geq 1$;
$$\| \exp (Nw) \| (\Omega_k ) \leq O(y^k) B .$$
Thus,
$$ \| \epsilon (w) \| = B \exp (- 2 \pi y) O(y^k) \quad \text{for some} \
k\leq n.$$
Similarly,
$$\align
\delta (w) &= \frac {\partial }{\partial w} ( \Omega (w) - \O (w) ) \\
   &= \frac {\partial }{\partial w} \sum_{k=1}^{\infty} \exp (2 \pi i k w)
 \exp (Nw) ( \Omega_k ) \\
&= \frac {\partial }{\partial w} \sum_{k=1}^{\infty} \exp (2 \pi i k w) ( 1
+ Niy + \dots + \frac { (N i y)^n }{n!} ) \Omega_k \\
&= \sum_{k=1}^{\infty}  2 \pi i k \exp (2 \pi i k w) ( 1+ Niy + \dots + \frac
 { (N i y)^n }{n!} ) \Omega_k  \\
&\quad  + \exp (2 \pi i k w) ( -N  - N^2 i y
 + \dots + \frac {N^n (i)^{n+1} y^{n-1} }{ (n-1)! } )\Omega_k  ) \\
&\quad +\exp (2 \pi i k w) ( 1+ Niy + \dots + \frac
 { (N i y)^n }{n!} )  (-N)\Omega_k
\endalign$$
Hence,
$$ \| \delta (w) \| \leq O(y^k) \exp ( - 2 \pi y ) C \quad\text {for some} \
  k \leq n \ \text {and} \  C>0 .$$
We have prove that both $\| \epsilon (w) \| $ and $\| \delta (w) \| $ have
exponential decay with respect to $y = Im (w)$.

Next, we look at the difference between the Weil-Petetsson arc length
integrands. Recall;
$$\align
R( \Omega (w) ) &= - \bar\partial w \partial w \log (\sqrt {-1})^n
Q ( \Omega (w), \bar\Omega
 (w) ) \quad \text {and} \\
R( \O (w) ) &= - \bar\partial w \partial w \log (\sqrt {-1})^n Q ( \O (w),
\bar\O (w) ).
\endalign $$
Note,from here on the form $Q( \  , \  ) $ includes $(\sqrt {\-1})^n $
So the difference of integrands is
$$ \sqrt{ R( \Omega (w)) } - \sqrt{ R( \O (w)) } = \frac {R( \Omega (w))
- R( \O (w))}{ \sqrt{R( \Omega (w)) } + \sqrt{ R( \O (w)) } } = (*)$$
The numerator become
$$\align
R( \Omega (w) ) - R( \O (w) )=&- \frac { Q(\Omega , \bar\Omega )
Q(\partial\Omega , \partial\bar\Omega )
- Q(\Omega , \partial\bar\Omega )Q( \partial\Omega , \bar\Omega ) }{
Q(\Omega , \bar\Omega )^2 } \\
&+ \frac {Q(\O ,\bar\O )Q(\partial\O , \partial
\bar\O ) - Q(\O ,\partial\bar\O )Q(\partial\O ,\bar\O )}
{ Q(\O ,\bar\O )^2 } = (**)
\endalign $$

where, $ \Omega = \Omega (w),\quad \O = \O (w) $.

By substituting $\Omega (w) = \O (w) - \epsilon (w) $ and
$\partial\Omega (w) = \partial\O (w) - \delta (w),$ we get;
$$\align
(**)=&- \frac {Q(\O +\epsilon  ,\bar\O + \bar\epsilon )Q(\partial\O +
\delta , \partial
\bar\O + \bar\delta ) - Q(\O + \epsilon  ,\partial\bar\O + \bar\delta
 )Q(\partial\O + \delta ,\bar\O + \bar\epsilon )}
{ Q(\O +\epsilon   ,\bar\O +\bar\epsilon )^2 } \\
&+ \frac {Q(\O ,\bar\O )Q(\partial\O , \partial
\bar\O ) - Q(\O ,\partial\bar\O )Q(\partial\O ,\bar\O )}
{ Q(\O ,\bar\O )^2 }
\endalign $$
After simplification we obtain;
$$
(**) = \frac 1{ Q(\O , \bar\O )^2 Q(\O +\epsilon , \bar\O + \bar\epsilon )^2 }
  F(w),$$
where each terms in $F( w) $ involves $\epsilon $, $\delta $, or both. Hence
$F(w) $ has a exponential decay.
Finally we simplify the entire $\sqrt {R(\Omega (w) )} - \sqrt {R( \O (w)) }$
and obtain;
$$
(*) = \frac { F(w) }{ (Q( \O , \bar\O  ))^3 \sqrt {
Q((\O , \bar\O )Q(\partial\O ,\partial\bar\O ) - Q(\O ,\partial\bar\O )
Q(\partial\O ,\O ) } + G(w) },
 $$
Where each term of G(w) involves $\epsilon $, $\delta $, or both.
Since $Q( \O, \bar\O )$ has the order of some power of $y$, the entire
term is dominated by exponential decay.

Thus we have showed that the difference between the two length is finite
along any vertical ray.

Now suppose the nilpotent orbit has finite length along given
parametric curve $ \alpha (t) : {\bold R} \longrightarrow \h $. Then the
difference between the
two length is given by;
$$\lim_{t \to \infty}\int_{t_0}^t (\sqrt {R(\Omega ( \alpha  (t) ) } - \sqrt{
R(\O (\alpha  (t)
)})  \sqrt  {(x'(t))^2 +(y'(t))^2 } dt \quad$$
$$\text{ where}\quad \alpha (t) = x(t) + i y(t)$$
By assumption the total variation of  $$ \sqrt  {(x'(t))^2 +(y'(t))^2}$$ has to
be finite.
Hence the difference is finite and the Weil-Petersson metric distance
of $\tilde\phi (w) $ has to be finite.

On the otherhand , if the nilpotent orbit has infinite length along
given parametric curve $\alpha (t) $, then the Weil-Petersson
metric length of the nilpotent orbit along any given
 parametric curve is infinite.
The difference of
two intengrand is,
$$ \sqrt{ R( \O (w)) } - \sqrt{ R( \Omega (w)) } = \frac {R( \O )
- R( \Omega)}{ \sqrt{R( \Omega) } + \sqrt{ R( \O ) } } \leq
\frac {R( \O )
}{ \sqrt{R( \Omega) } + \sqrt{ R( \O ) } }.
$$
Since $\sqrt{R( \Omega )}$ is dominated by $\sqrt{R (\O )}$,
the difference is dominated by a half of the nilpotent orbit
distance integrand.
 Hence  the distance along $\tilde\phi (w)$ has to be infinite.

Thus we complete the proof.

\Bigskip
\head 2. Limiting  mixed Hodge structure \endhead
\bigskip

\noi Our next step is to determine for which limiting mixed Hodge structure
 the central fibre is at
finite distance in the  Weil Petersson metric.
We have to define some terminology first. These
definitions are found in Chapter 4 and Chapter 5 of \cite{GrTA}.
\smallskip
\noi {\bf Definition 2.1} Given a variation of Hodge structure of
 weight $n$ over
$\D^*$, and a nilpotent monodromy $N$ such that $N^{n+1}=0$, there exists a
unique ascending filtration $\{ W_l \}$ of $H_{\bold Q}$, called the
{\it monodromy filtration},
$$0\subset W_0 \subset W_1 \subset \ldots \subset W_{2n}=H_{\bold Q},$$
\noi satisfying
$$ N\: W_i \rightarrow W_{i-2} $$
$$ N^k \: W_{n+k}/W_{n+k-1} \cong  W_{n-k}/W_{n-k-1}.$$
\smallskip
\noi Note that the monodromy filtration is uniquely determined by these
properties.

\noi {\it Caution }  All filtrations are on the complexification
of the primitive rational cohomology in dimension $n$. We are about
to introduce a separate notion of primitivity related to the weight
filtration. Any further reference to primitivity refers to
 newly introduced primitivity.
\smallskip
\noi {\bf Definition 2.2}  Given a nilpotent element $N \in \g_{\bold
R}$, such that $N^{n+1}=0$  its
 monodromy weight filtration $W=W(N)$,
 and each graded piece $Gr_l^W=\frac {W_l}{W_{l-1}}$,
the {\it primitive subspaces}  $P_{n+j}\subset Gr_{n+j}^W$ are  defined as
follows:
$$P_{n+j}={\bold ker} \{N^{j+1}\: Gr_{n+j}^W\longrightarrow Gr_{n-j-2}^W
\}\quad\text {if\ }
\quad j\geq 0$$
$$P_{n+j}={0}\quad\text {if\ }\quad j<0.$$
\smallskip
We have a Lefschetz decomposition:
$$Gr_l^W=\underset {j\geq0}\to\oplus N^j(P_{l+2j})$$
Given $n$,  $P_{n-1}=P_{n-2}= \ldots =P_1=P_0=0$, and

Next we define the mixed Hodge structure.
\smallskip
\noi {\bf Definition 2.3} {\it A mixed Hodge structure} $\{ H_{\bold Q}, F^p ,
W_l \} $ of weight $n$ consists of an ascending weight filtration, defined
over {\bf Q},
$$0 \subset W_0 \subset W_1 \subset W_2 \subset W_3 \subset \ldots \subset
W_{2n}= H_{\bold C},$$

\noi such that the Hodge filtration induces the pure Hodge structure of weight
$m$ on the graded piece $Gr_m =W_m/W_{m-1} $ of  the weight filtration for
each $m=0,\ldots ,2n$. The induced Hodge filtration on $Gr_m$ is
$$F^p(Gr_m)=(F^p \cap W_m)/(F^p \cap W_{m-1}) .$$
\smallskip
\noi We have  polarized versions of these definitions.
\smallskip
\noi {\bf Definition 2.4  } {\it
A polarized mixed Hodge structure} is a triple $(W,F,N)$, where $W$ is an
increasing filtration of $H$, $F\in\check D$ and
$N \in \g_{\bold R} $ is a rational  nilpotent element such that $N^{n+1}=0$,
satisfying:
\roster
\item $W$ is the monodromy weight filtration of $N$.
\item $(W,F)$ is a mixed Hodge structure; that is, for every $l\geq0$ the
filtration induced by $F$ on $Gr_l^W$ is a
      Hodge structure of weight $l$.
\item $N(F^p) \subset F^{p-1}:\quad\text { for all}\quad 0 \leq p \leq n.$
\item For $j\geq0$, the Hodge structure induced by $F$ on the primitive
subspace $P_{n+j}$ is polarized by the bilinear
   form $Q_j=Q(.,N^j.)$.
\endroster
\smallskip

We represent the limiting mixed Hodge structure of $H^n$ pictorially as
follows.

$$H^{n,n}$$
$$H^{n,n-1},H^{n-1,n}$$
$$....................................................$$
$$H^{n,0},H^{n-1,1},H^{n-2,2},\ldots,H^{2,n-2},H^{1,n-1},H^{0,n}$$
$$....................................................$$
$$H^{2,0},H^{1,1},H^{0,2}$$
$$H^{1,0},H^{0,1}$$
$$H^{0,0}$$
Each row corresponds to the Hodge structure of each graded piece.
Furthermore,
 we have isomorphisms
$$N:Gr_{n+1} \cong Gr_{n-1},\quad N^2:Gr_{n+2} \cong Gr_{n-2},\quad
N^3:Gr_{n+3} \cong Gr_{n-3},\ldots ,$$
$$N^{n-1} :Gr_{2n-1} \cong Gr_{1},\quad N^n :Gr_{2n} \cong Gr_{0} .$$
\smallskip
Hence, there is a horizontal symmetry (by Hodge structure) and a
vertical symmetry (by isomorphisms $N^i, $ for $1 \leq i  \leq n $ ) in this
pictorial
representation of
the limiting mixed Hodge structure.

Since  $F^n_{\infty} $  has rank $1$, exactly one of the graded pieces
($ Gr_{n+k}, 0 \geq k \geq n$) is nonzero.
 There are  $n+1$ cases.

\noi Due to the following theorem of Schmid, the limiting filtration does
induce
the polarized limiting mixed Hodge structure.

\smallskip
\noi {\bf Theorem (Schmid)  } {\it
Let $\Cal O(w)=e^{wN}F_{\infty}$ be a one variable nilpotent orbit where
$F_{\infty}$ is a limit filtration. Then
$(W(N),F_{\infty},N)$ is a polarized mixed Hodge structure.}
\smallskip
\noi {\it Proof  }    See  theorem 6.16 of \cite{Sch}
\smallskip

Now we can compute the polarized mixed Hodge structure associated to the limit
filtration and impose the condition
of the finite distance case. Since our variation of Hodge structures
 comes from the degeneration space of a family of Calabi-Yau n-folds, the
dimension
of $F_{\infty}^n$ is $1$.
\smallskip
\noi {\bf Theorem 2.5  } {\it
Only  one  case  satisfies the  finite distance condition on $F_{\infty}^n$;
 i.e.,
$$Q(\alpha,N^i\bar\alpha)=0,\,\text {for all\ }\,i>0$$
$F^n_{\infty}$ is at finite distance if and only if $F^nGr_n(H^n(\Xs)
=H^{n,0}(\Xs)$.}
\smallskip
\noi {\it Proof  }
Since $H^{n,0}(\Xs)$ has rank $1$,
 $F^n_{\infty}=F^nGr_{n+k}(H^n(\Xs))$,is different from $0$ for
precisely one $k$, where $k,\ 1\geq
k \geq n$. Then,
$$F^n_{\infty}\subset P_{n+k},\quad\text {and}\quad Q(F^n_{\infty},
N^k\bar F^n_{\infty})>0,
$$
by the definition of polarized mixed Hodge structure, where
$$F^n_{\infty}=F^nGr_{n+k}(H^n(\Xs))$$
is polarized by
$Q(.,N^k.)$ which satisfies {\bf HR}-bilinear relations. Thus
$$N^{k+1}F^n_{\infty}=0.$$
Hence $F^n_{\infty}=F^nGr_{n+k}$ is at finite distance if and only if
$k=0$, i.e.,
$$F^nGr_n(H^n(\Xs)=H^{n,0}(\Xs)\neq 0.$$
\line{\hfill \qed}
\smallskip
\noi {\bf Corollary 2.6  } {\it
The pictorial presentation of the limiting mixed Hodge structure
at finite distance is the following. (See right after the proof.)

Furthermore, $N^{n-1}=0$ for this limiting mixed Hodge structure.}
\smallskip
\noi {\it Proof  }
The picture is clearly determined once we set the non-zero part of
$F^n_{\infty}$ to be $F^nGr_n(H^n(\Xs))=H^{n,0}(\Xs)\ $.
  From the picture , each graded piece becomes
$$Gr_{2n}=P_{2n}=0,Gr_{2n-1}=P_{2n-1}=0$$
Hence the first  non-zero primitive piece is
$$\split
P_{2n-2}
  &={\bold ker} \{ N^{n-2+1} \: Gr_{n+n-2} \longrightarrow Gr_{n-(n-2)-2}\} \\
  &={\bold ker} \{ N^{n-1} \: Gr_{2n-2} \longrightarrow Gr_0 \}
\endsplit
$$
Hence, by the Lefschetz   decomposition, $N^{n-1}$ annihilates everything.

\line{\hfill \qed}
$$0$$
$$0,0$$
$$0,X,0$$
$$0,Y,Y,0$$
$$...............$$
$$1,Z,..................Z,1$$
$$...............$$
$$0,Y,Y,0$$
$$0,X,0$$
$$0,0$$
$$0$$

\bigskip
\head 3. Central Fibre and Clemens-Schmid exact sequence \endhead
\bigskip

\noi In order to interpret the finite distance condition of a Calabi-Yau
degeneration we need to look at the central fibre. In fact, every variety
carries a canonical functorial mixed Hodge structure; as proved
in a more general setting by Deligne \cite{DeH2} and \cite{DeH3}. Before
looking
at the mixed Hodge structure of a central fibre, we need to understand
its cohomology.
\smallskip
\noi {\it Remark 3.1  }  In the case of a semistable degeneration, one
can construct a retraction from the total space to the central fibre; i.e.,
$r \: \XT \longrightarrow \X0 $, which induces  isomorphisms
$$r^* \: H^m(\X0,{\bold Q}) \longrightarrow H^m(\XT, {\bold Q}),\quad
r_* \: H_m(\XT,{\bold Q}) \longrightarrow H_m(\X0, {\bold Q})$$
This construction is due to Clemens in \cite{Cl}.
The cohomology of a central fibre is described by a Mayer-Vietoris
type spectral sequence. We follow the construction given in
Chapter 6 of \cite{GrTA}.

\smallskip
\noi {\it Construction 3.2  }  Let $\X0= \sum X_i$ where each $X_i$ is
a  smooth irreducible component of $\X0 $ and all components
intersect  transversely each
other.
 Let $X_{i_0 \ldots i_p}$ denote $X_{i_0} \cap X_{i_1} \cap \ldots \cap
X_{i_p} $. We define the {\bf codimension $p$ stratum } of $\X0$ as
$$ \XT^{[p]} = \underset {i_0 < \ldots < i_p} \to\bigsqcup  X_{i_0 \ldots i_p}
\quad\text {(disjoint union)}.$$ There is a natural inclusion map $\iota_p \:
\XT^{[p]}\longrightarrow \XT.$

By \cite{Cl} we can choose an open cover $\Cal U$ of a neighborhood of $\X0$ in
$\XT$ satisfying the following condition;

\roster
\item For each $U \in \Cal U$, the degeneration map $\pi$ restricted to
$U$ is given by $s=x_1 x_2 x_3 \ldots x_k$ in suitable local cordinates.
\item $\check H^*(\Cal U\cap\X0 , {\bold Q}) \cong \check H^*(\X0 , {\bold Q})$
\item $\check H^*(\iota_p^{-1} (\Cal U),{\bold Q}) \cong \check
H^*(\XT^{[p]},{\bold Q})$.
\endroster

\noi Then we can define $E_0$ terms and differentials as followings;
$$E_0^{p,q}=\check C^q (\iota_p^{-1} (U),{\bold Q})$$
$$\align
&d \: E_0^{p,q} \rightarrow E_0^{p,q+1}, \qquad \text {the Cech coboundary,} \\
&\delta \: E_0^{p,q} \rightarrow E_0^{p+1,q} , \qquad \text {the combinatorial
coboundary induced by,}\\
&\delta\phi ( U\cap X_{j_0\ldots j_{p+1}}) = \underset\alpha\to\sum
(-1)^{\alpha}\phi (U\cap X_{j_0\ldots \hat j_{\alpha} \ldots j_{p+1}}).
\endalign$$

\noi This defines a bigraded complex and we get an associated spectral
sequence.

We have the following result of  the
convergence of this spectral sequence  by  Griffiths and Schmid in \cite
{GrSc}.
\smallskip
\noi {\bf Theorem    } {\it
  The spectral sequence given above degenerates at $E_2$, and  converges
to $\check H^*(\X0, {\bold Q})$.}
\smallskip
\noi {\it Proof  }  See \cite{GrSc}.
\smallskip
\noi Theorem comes from the de Rham analogue of our spectral sequence.
Let
$$A^{p,q} := A^q(\XT^{[p]})\quad \text {be the  complex  of}\quad
C^{\infty}\quad
 q-\text
 {forms on } \XT^{[p]},$$
with two differentials,
$$\align
&d \:A^{p,q} \rightarrow A^{p,q+1}, \qquad \text {the exterior differential,}
\\
&\delta \: A^{p,q} \rightarrow A^{p+1,q} , \qquad \text {the combinatorial
coboundary induced by,}\\
&\delta\phi_{j_0\ldots j_{p+1}} = \underset\alpha\to\sum
(-1)^{\alpha}\phi_{j_0\ldots \hat j_{\alpha} \ldots j_{p+1}}|_{X_{j_0\ldots
j_{p+1}}}.
\endalign$$
Then this gives a bigraded complex and has an associated spectral
sequence. It converges to $H^*_{\bold  DR} (A^{*,*})$
 and degenerates at $E_2$.  By the de Rham Theorem for an algebraic
variety with normal crossings, we have
$$H^*_{\bold DR} (A^{*,*}) \cong H^*(\X0, {\bold C}).$$
At the $E_1$ level, we have a commutative diagram:
$$\diagram
E_r^{p,q}\rto^{d_r}\dto_i  &  E_r^{p,q+1}\dto _i \\
E_r^{p,q}\otimes{\bold C}\rto^{d_r}\dto_{\cong} &  E_r^{p,q+1}\otimes{\bold C}
\dto_{\cong} \\
\ ^{\bold DR}E_r^{p,q}\rto^{d_r}  & \ ^{\bold DR}E_r^{p,q+1} \\
\enddiagram
$$
This ensures that the topological spectral sequence degenerates at
$E_2$.

Now we are ready to define the mixed Hodge structure on the cohomology
of the central fibre by means of a spectral sequence.
First, this spectral sequence
gives a means to put a {\bf weight
filtration} on $H^{*}(\X0, {\bold Q})$ as follows:

\noi define
$$ W_k =\bigoplus_{q \leq k} E_0^{*,q} ,$$
and let $W_k( H^m)$ be the induced filtration. Then this will give an
increasing filtration of length $m$ on $H^m$, and the graded pieces
will be
$$Gr_k(H^m)=E_2^{m-k,k} \quad\text {if }\quad 0\leq k \leq m\quad \text
 {and}$$
$$Gr_k(H^m)=0 \quad
\text {if}\qquad k<0 \quad\text {or}\quad k>m.$$

Next we put a filtration on the complexification of the original
spectral sequence which is isomorphic to  the usual de Rham sequence, i.e.,
at $E_0$ level,
$$ F^p(A^{k,l})=\bigoplus_{r\geq p} A^{r,l-r}(\XT^{[k]})$$
This induces filtrations on each $E_r$ term, in particular, at the $E_1$
level we get the   Hodge filtration:
$$\ ^{\bold DR}E_1^{k,l}=H^l_{{\bold DR}}(\XT^{[k]})$$
Furthermore, $d_1$ is  a morphism of Hodge structures. So, together
with the weight filtration, $H^*(\X0, {\bold C})$ carries
a mixed Hodge structure i.e., it carries a Hodge structure of weight $k$
at each $G_k$ piece.

To illustrate the filtrations we  represent the $E_1$ level of
the spectral sequence associated to an $n-$dimensional central
fibre.

$$\align
E_1^{n,0} &=H^{n,0}(\XT^{[0]})\oplus H^{n-1,1}(\XT^{[0]})\oplus\ldots\oplus
H^{1,,n-1}(\XT^{[0]})\oplus H^{0,n}(\XT^{[0]})\\
E_1^{n-1,1} &=H^{n-1,0}(\XT^{[1]})\oplus H^{n-2,1}(\XT^{[1]})\oplus\ldots\oplus
H^{1,n-2}(\XT^{[1]})\oplus H^{0,n-1}(\XT^{[1]})\\
E_1^{n-2,2} &=H^{n-2,0}(\XT^{[2]})\oplus H^{n-3,1}(\XT^{[2]})\oplus\ldots\oplus
H^{1,n-3}(\XT^{[2]})\oplus H^{0,n-2}(\XT^{[2]})\oplus
\endalign$$
$$.............................................................................
..$$
$$\align
E_1^{1,n-1} &=H^{1,0}(\XT^{[n-1]})\oplus H^{0,1}(\XT^{[n-1]})\\
E_1^{0,n} &=H^{0,0}(\XT^{[n]})
\endalign$$
The $E_2$ terms  give the lower half triangle of  the picture
of the mixed
Hodge structure  we used in Section $2$. In other words,
the mixed Hodge structure of the central fibre has all $0$'s for
the upper half.

For our purpose, we need to look at  $F^n$, the $n$th filtration of weight
$n$ Hodge filtration. Since $d_1$ respects the Hodge filtration , we
get $E_2^{n,0}$ by looking at
$$0\longrightarrow F^nE_1^{n,0}\longrightarrow F^nE_1^{n,1}\longrightarrow 0,$$
where
$$\align
&F^nE_1^{n,0} \cong H^{n,0}(\XT^{[0]})\\
          &F^nE_1^{n,1} \cong H^{n,1}(\XT^{[1]}).
\endalign$$
Notice, since $\XT^{[1]}$ is a first level stratum, its complex dimension
is $n-1$, i.e., it can have  at most $n-1$ non-trivial filtration
depth. Hence,
$$F^nE_1^{n,1}=0\qquad\text {and}\qquad F^nE_2^{n,0}=F^nE_1^{n,0}$$
To interpret the finite distance condition, this is the crucial component
of the mixed Hodge structure of a central fibre.

 The final
key construction to connect the mixed Hodge structure induced by
the limiting filtration to the canonical mixed Hodge structure of a central
fibre is the  Clemens-Schmid exact sequence.
The most elaborate description can be found in \cite{Cl}, and
we again more or less follow the one given in \cite{GrTA} Ch6.
Since maps of the Clemens-Schmid exact sequence are not just maps
of homology and cohomology groupes but morphisms of mixed Hodge
structures of appropriate types,  we need first to define those
morphisms precisely.
\smallskip
\noi {\bf Definition 3.3  }  A {\it weighted Q vector space} is a
{\bf Q} vector space $H$ together with an increasing filtration of
{\bf Q} subspaces
$$0\subset\ldots\subset W_k(H)\subset W_{k+1}(H)\subset\ldots\subset H. $$
Then a {\rm morphism of weighted vector spaces of type (r,r)} is a
linear map $\phi \: H \rightarrow H' $ such that
$$\phi( W_k(H))= W_{k+2r}(H') \cap {\bold Im}\phi.$$
\smallskip
Both the monodromy weight filtration given in section 2 and the weight
filtration given in this section  are examples of weighted {\bf Q}
vector spaces. The monodromy operator $N$ we have been looking
at is an example of a type $(-1,-1)$ morphism from $H^m(\Xs ,{\bold Q})$
 to itself.

We also need to recall the maps appearing in the Clemens-Schmid
exact sequence.
\roster
\item Let $i \: \Xs\hookrightarrow\XT$  be a natural inclusion of the
generic fiber into the total space. Since we know $H^m(\XT) \cong H^m(\X0)$,
this induces
$$i^* \: H^m(\X0) \rightarrow H^m(\Xs).$$
\item Define $\alpha \: H_{2n+2-m}(\X0) \rightarrow H^m(\X0)$ to be
the composite

$$H_{2n+2-m}(\XT) \overset {P.D}\to\rightarrow H^m(\XT, \partial\XT)\rightarrow
H^m(\XT)$$
where $P.D$ is the Poincar\'e duality map. Again, we  identify
cohomology of the central fibre with that  of the total space.
\item Define $\beta \: H^m(\Xs) \rightarrow H_{2n-m}(\X0)$ to be the
 composite
$$H^m(\Xs) \overset {P.D}\to\rightarrow H_{2n-m}(\Xs) \overset
 {i_*}\to\rightarrow H_{2n-m}(\Xs)$$.
\endroster
For further details see \cite{Cl}.
\smallskip
\noi {\bf Theorem of Clemens-Schmid exact sequence  }
 {\it  The maps $\alpha , i^* , N , \beta $ are morphisms of
 weighted vector spaces of type $(n+1,n+1)$, $(0,0)$, $(-1,-1)$,
 $(-n,-n) $, respectively,
and the sequence
$$\align
\rightarrow H_{2n+2-m}(\X0) \overset\alpha\to\rightarrow H^m(\X0)
\overset {i^*}\to\rightarrow H^m(\Xs) &\overset N\to\rightarrow H^m(\Xs)
\overset\beta\to\rightarrow\\
&\rightarrow H_{2n-m} \rightarrow^{\alpha} H^{m+2} \rightarrow \ldots
\endalign$$
is exact.  Furthermore the morphisms
in the Clemens-Schmid sequence are morphisms of  mixed Hodge structures
of the appropriate types.}
\smallskip
\noi {\it Proof and construction  } {\it See \cite{Cl}.}
\smallskip
\noi For our purpose, we need to look at the exact sequence around $m=n$, i.e.,
$$\align
\rightarrow H_{2n+2-n}(\X0) \overset\alpha\to\rightarrow H^n(\X0)
\overset {i^*}\to\rightarrow H^n(\Xs ) &\overset N\to\rightarrow H^n(\Xs)
\overset\beta\to
\rightarrow\\
&\rightarrow H_{2n-n} \overset\alpha\to\rightarrow H^{n+2} \rightarrow \ldots
\endalign$$
Furthermore, we can  restrict the Clemens-Schmid exact sequence to the
relevant filtered part of the  graded pieces.

\noi Here we need one more  definition;  we have to
define the weight filtration on homology that gives the right definition
of the graded pieces of homology.  That is,
$$W_{-k}(H_m)=Ann(W_{k-1}(H^m))= \{ h\in H_m | (W_{k-1}(H^m), h)=0 \}.$$
\noi Then,
$$Gr_k(H_m) \cong (Gr_{-k}(H^m))^* \qquad \text {so that}$$
$$\qquad Gr_k(H_m)=0\quad\text {if}\quad k<-m\quad \text {or}\quad k>0.$$

In our case, $H^n(\X0)$ carries the mixed Hodge structure we described in
this section;  that is, the canonical one of the central fibre, and $H^n(\Xs)$
carries the mixed Hodge structure in section 2; that is, the limiting mixed
Hodge structure.

\noi Hence,

\roster
\item  $i^* \: H^n(\X0)\rightarrow H^n(\Xs)$ sends each piece of the mixed
Hodge
structure of a central fibre to the lower half of the limiting mixed Hodge
structure without any filtration index shift; i.e.,
$$F^p(Gr_m H^n(\X0))\rightarrow  F^p(Gr_m H^n(\Xs))\qquad \text {where}
\qquad  0\geq p,\quad m \geq n .$$
\item  $N \: H^n(\Xs)\rightarrow H^n(\Xs) $ sends each piece of the limiting
mixed Hodge Structure into itself but shifts  the Hodge filtration
index by $(-1,-1)$ and monodromy weight filtration index by $-1-1=-2$
; i.e.,
$$F^p(Gr_m H^n(\Xs))\rightarrow F^{p-1}(Gr_{m-2}H^n(\Xs)).$$
\endroster

Now we can restrict the Clemens-Schmid exact sequence to the following graded
pieces.
$$\align
F^{-1}Gr_{n-(2n+2)}H_{(2n+2)-n}(\X0)\overset\alpha\to\rightarrow
F^nGr_n H^n(\X0)&\overset {i^*}\to\rightarrow\\
\rightarrow F^nGr_n H^n(\Xs)&\overset N\to\rightarrow
F^{n-1}Gr_{n-2} H^n(\Xs)\overset\beta\to\rightarrow
\endalign$$
and
$$ F^{n-1}Gr_{n-2} H^n(\Xs)=0$$
since the  $n-2$ graded piece carries only
up to $n-2$ depth Hodge filtration and there will be no $n-1$ deep
filtration piece.

\noi Furthermore, by a linear algebra fact about mixed Hodge
 structures (see \cite{GrSc}),
$$F^{-1}Gr_{-n-2}H_{n+2}
(\Xs)=Ann(F^2Gr_{n+2} H^{n+2}(\Xs)),$$
 and it is $0$
 since we have only up to $n-$forms for K\"ahler manifolds.

 Hence, we have
$$F^nGr_nH^n(\Xs) \cong F^nGr_nH^n(\X0).$$

\noi Now we are ready to state and prove the main result.
\smallskip
\noi {\bf Theorm 3.4 }  {\it  A semistable degeneration of a Calabi-Yau
$n$ dimensional compact complex variety is at  finite distance with
respect to the Weil-Petersson metric if and only if $H^{n,0}$ of one
of the components of the central fibre has positive rank.}
\smallskip
\noi {\it Proof  }  Recall the result from this section, where we computed
the $E_1$ terms of the spectral sequence of the central fibre.
$$\align
 H^{n,0} (\XT^{[0]}) = F^nE_1^{n,0}\cong F^nE_2^{n,0}&\overset
 {(*)}\to = F^nGr_nH^n(\X0)\\
&\text {where}\quad (*)\quad\text {is by
definition.}
\endalign$$
On the other hand,  a degeneration is at finite distance if and only if
$$F^nGr_nH^n(\Xs) \neq 0 , \quad\text {i.e.,}$$
$$F^nGr_nH^n(\Xs) \cong F^nGr_nH^n(\X0) \cong H^{n,0} (\XT^{[0]})\neq 0.$$
By definition, $\XT^{[0]}$  corresponds to the disjoint union of
irreducible pieces of the central fibre of the degeneration.
 Hence,  exactly one
of the components of the central fibre has  $h^{n,0}\neq 0$.

\line{\hfill \qed}
\smallskip
\noi {\it Remark } {\it Since the dimension of $F^n_{\infty}$ is exactly
one, actually $h^{n,0}=1$.}
\Bigskip
\head 4. Examples \endhead
\bigskip

\noi Originally we were looking at only the Calabi-Yau $3-$fold case, which is
 the main ingredient in the phenomena called `mirror symmetry' in mathematical
physics. Since the result applies in any dimension, we first check this
with the known result of the low dimensional degeneration problem.
\smallskip
\noi {\it Case 1.  }   Degeneration of Curves
\smallskip
\noi Complex curves which have trivial canonical bundles are complex
tori. It is  a classical result that the genus of the central fibre
of a semi-stable degeneration is strictly less than the genus of
the generic fibre. Hence,  each irreducible component of the central fibre
has to consist of ${\bold CP^1}$'s and none of $H^{1,0}$ is non-zero.
Hence the degeneration has to be at  infinite distance, which coincides with
the known result.
\smallskip
\noi {\it Case 2.  }   Degeneration of Surfaces
\smallskip
\noi Smooth complex surfaces with trivial canonical
bundles are K3 surfaces. The semi-stable degeneration of K3
surfaces is classified by the following theorem.
\smallskip
\noi {\bf Theorem  (Kulikov \cite{Kul}, Persson and Pinkham \cite{PP})}
{\it A semi-stable degeneration of K3 surfaces is birational to one
for which the central fibre $\X0$ is one of the following three;
\roster
\item $\X0$ is smooth K3 surface.
\item $\X0 =X_0 \cup X_1 \cup \ldots \cup X_{k+1}$. $X_i$ meets only
$X_{i\pm 1}$, and each $X_i \cap X_{i+1} $ is an elliptic curve.
$X_0$ and $X_{k+1}$ are rational surfaces, and for $1 \leq i \leq k$,
$X_i$ is ruled with $X_i \cap X_{i+1}$ and $X_i \cap X_{i-1}$ sections
of the ruling.
\item All components of $\X0$ are rational surfaces, $X_i \cap (
\underset {j \neq i}\to\cup X_j )$ is a cycle of rational curves, and
the dual graph of the configuration is homeomorphic to $S^2$.
\endroster}
\smallskip
\noi  Recall (cf. \cite{GrHr} ) both rational and elliptic ruled
surfaces have {\bf Kodaira number} $-1$, i.e.,
$$ 0= h^0(X_i, \Cal O(K_{X_i})) = h^{2,0}(X_i) .$$
Hence only smooth central fibre gives  finite distance.
\smallskip
\noi {\it Case 3.  }   Nodal degeneration of Calabi-Yau $n$-folds,
 for $n \geq 3$
\smallskip
To get a semi-stable degeneration we have to blow-up all nodes,
so that we will get a central fibre with  normal crossings.
Let $\pi\:\XT\longrightarrow\ \D\ $ be a degeneration of Calabi-Yau
$n$-folds.
Suppose that the special fibre $\X0$
is a $n-$dimensional compact complex variety with $m$ nodes. Then the
full blow-up of $m$ nodes, denoted $\X0'$, has the following configuration.

There are $m$  exceptional divisors of the total space $\XT,$ ${\bold C}P^n $,
and a proper
transform of $\X0$, $\bar\X0$.
Each ${\bold C}P^n$ intersects with $\bar\X0$ once and only once;
none of ${\bold C}P^n$
intersects each other. There are no self-intersections.

\smallskip
\noi {\bf Proposition 4.1  }  {\it The blow-up of nodes are at finite
distance with respect to Weil-Petersson metric. Hence the special
fibre with nodes is  at finite distance.}
\smallskip
\noi {\it Proof  } Let $\pi'\:\XT'\longrightarrow\ \D\ $be the
degeneration of the blow-up of nodes $\XT'$, which is clearly semi-stable.
If we denote $f\:\XT'\longrightarrow\ \XT\ $ to be the full blow-up
map, then $\pi' = \pi\circ f$.

Since each generic fibre is a Calabi-Yau fibre,
the canonical divisor of $\XT$ is
concentrated on the blow-up of the special fibre $\X0'$. Hence we have
$$K_{\XT'}=kL_{\bar\X0} + \Sigma k_i L_{D_i} \tag1 $$
where $L_{\bar\X0}$ is the line bundle corresponding to $\bar\X0$ as an
effective divisor in $\XT'$ and $L_{D_i}$ is the line bundle corresponding
to each exceptional divisor $Di \cong {\bold C}P^n $ in $\XT'$.

Let $H_i$ denote the hyper-plane bundle over $D_i$. Then
$$\align
L_{D_i}|_{D_i} &= N_{D_i}^* = - H_i \in {\bold Pic}(D_i), \tag2 \\
K_{D_i} &= K_{{\bold C}P^n} = -(n+1)H_i . \tag3
\endalign$$

\noi {\it Claim }
$$ L_{\bar\X0} + 2\Sigma L_{D_i} = 0.$$
\noi {\it Proof of Claim}  Since $\XT'$ is the full blow-up of nodes, each
exceptional
divisor has multiplicity $2$, i.e.,
$$\bar\X0 + 2\Sigma D_i = f^{-1}(\X0).$$
If we look at the trivial bundle of $\XT'$, it has a section given by
$$ x'\longmapsto (x', \pi'(x') ) = (x', \pi\circ f (x')).$$
Notice the zero locus of this section is exactly $\pi'(f^{-1}(\X0))$.
Hence the corresponding line bundle is trivial.

\line{\hfill \qed}

Hence we have
$$L_{\bar\X0} = -2\Sigma L_{D_i}.$$
We also get other equalities of line bundles:
$$L_{\bar\X0}|_{D_i} = -2\underset j\to\Sigma L_{D_j}|D_i =-2L_{D_i}|D_i
=2H_i \tag4 $$
$$\text {since}\quad L_{D_j}|_{D_i} = L_{D_j \cap D_i} = L_{\emptyset} =0,$$
and
$$L_{\bar\X0}|_{\bar\X0} = -2\Sigma L_{D_i}|_{\bar\X0}= -2\Sigma L_{E_i}
\in {\bold Pic}(\bar\X0) tag5 ,$$
$$ \text {where}\quad E_i = \bar\X0 \cap D_i.$$
We are now ready to compute the canonical bundle of $\bar\X0$.

\noi By the adjunction formula for $D_i$ we have
$$K_{D_i} = K_{\XT'}|_{D_i} + L_{D_i}|_{D_i}.$$
Notice
$$\align
 K_{\XT'}|_{D_i} &= kL_{\bar\X0}|_{D_i} + \underset j\to\Sigma k_j
L_{D_j}|_{D_i} \qquad \text {by}\quad (1) \\
&= 2kH_i - k_i H_i \qquad \text {by}\quad (2),(4).
\endalign$$
Then we get
$$\align
-(n+1)H_i &= 2kH_i - k_iH_i - H_i \qquad \text {by}\quad (3) \\
&= (2k-k_i -1)H_i .
\endalign$$
Since the first Chern class gives the isomorphism:
$${\bold c}_1 \: {\bold Pic}({\bold C}P^n) \longrightarrow H^2({\bold C}P^n,
{\bold Z}),$$
We get
$$-(n+1)= 2k-k_i -1,\qquad \text {i.e.,}\qquad k_i= 2k+n.$$
Next, by the adjunciton formula for $\bar\X0$,
$$K_{\bar\X0} = K_{\XT'}|_{\bar\X0} + L_{\bar\X0}|_{\bar\X0}.$$
Then by $(1)$ and $L_{\bar\X0} = -2\Sigma L_{D_i}$, we get
$$\align
K_{\XT'} &= \Sigma k_i L_{D_i} -2k\Sigma L_{D_i} \\
 &=\Sigma (k_i - 2k)L_{D_i} \\
 &=\Sigma (2k+n-2k)L_{D_i} \\
 &=\Sigma n L_{D_i}.
\endalign$$
Hence
$$\align
 K_{\bar\X0} &= \Sigma n L_{D_i}|_{\bar\X0} + L_{\bar\X0}|_{\bar\X0} \\
   &= n\Sigma L_{E_i} -2\Sigma L_{E_i} \\
   &= (n-2)\Sigma L_{E_i}.
\endalign$$
Recall $E_i = \bar\X0 \cap D_i$. So $E_i$ is an effective divisor in $\bar\X0$.
Hence for $n \geq 3 $, $K_{\bar\X0}$ has a non-trivial global holomorphic
section,
 i.e., $H^{n,0} (\bar\X0) \neq 0$. By Theorem {\bf 3.4}  result follows.

\line{\hfill \qed}
\smallskip
\noindent {\it Remark }  This result coincides
with the lower dimensional cases.
For surfaces, the special fibre has to be smooth. Therefore, it is also
a K3 surface which has a  trivial canonical bundle.  If  the
fibres are curves, the
canonical bundle cannot have a holomorphic global section.

\vfill\eject
\centerline{{\bf REFERENCES}}
\bigskip
\widestnumber\key{CaGrTH}
\ref \key Ks \by Norihito Koiso
\paper Einstein metrics and complex structures
\jour Invent.-Math
\yr 1983 \vol 73 \pages 71--106
\endref
\ref \key CGH \by P. Candelas, P. Green and T. H\"ubsch
\paper Rolling among Calabi-Yau Vacua
\jour Nuclear Physics
\yr 1990 \vol B330 \pages 49--102
\endref
\ref \key Cl \by C. H. Clemens
\paper Degeneration of K\"ahler manifolds
\jour Duke Math J.
\yr 1977 \vol 44 \pages 215--290
\endref
\ref \key DeH2 \by P. Deligne
\paper Th\'eorie de Hodge II
\jour Pub. Math. IHES \yr 1972 \vol 40 \pages 5--57 \endref
\ref \key DeH3 \by P. Deligne
\paper Th\'eorie de Hodge III
\jour Pub. Math. IHES  \yr 1973 \vol 44 \pages 5--77 \endref
\ref \key Gr1 \by P. Griffiths
\paper Period of integrals on algebraic manifolds I, II
\jour Amer. J. Math \yr 1970 \vol 90  \pages 568--626
\moreref \pages 805--865
\endref
\ref \key GrCo \by M. Cornalba and P. Griffiths
\paper Some transcendental aspects of algebraic geometry
\inbook in Proceedings of Symposia in Pure Mathematics
\publ A.M.S. \yr 1975 \vol 29 \pages 3--110
\endref
\ref \key GrHr \by P. Griffiths and J. Harris
\book Principles of Algebraic Geometry
\publ A Wiley-Interscience Publication \yr 1978
\endref
\ref \key GrSc \by P. Griffiths and W. Schmid
\paper Recent developments in Hodge Theory: a discussion of
techniques and results
\inbook Discrete Subgroups of Lie Groups and Application to Moduli
\publ Oxford University Press \yr 1973 \pages 31--127
\endref
\ref \key GrTA \by P. Griffiths ed.
\book Topics in Transcendental Algebraic Geometry
\bookinfo Annals of Mathematical Studies
\publ Princeton University Press \yr 1984 \vol 106
\endref
\ref \key KS \by K. Kodaira
\book Complex Manifolds and Deformation of Complex structure
\publ Springer -Verlag \publaddr New York \yr 1986
\endref
\ref \key Kul \by V. Kulikov
\paper Degeneration of K3 surfaces and Enriques surfaces
\jour Math  the USSR Izvestija \yr 1977 \vol 11 \pages 957--989
\endref
\ref \key LA \by A. Landman
\paper On the Picard-Lefschetz transformations
\jour Trans A.M.S. \yr 1973 \vol 181 \pages 89--126
\endref
\ref \key Mum \by G. Kempf et al
\book Troidal Embeddings. I
\bookinfo Lecture Notes in Math., \yr 1973 \vol 339
\publ Springer -Verlag \publaddr Berlin and New York
\endref
\ref \key PP \by U. Persson and H. Pinkham
\paper Degeneration of surfaces with trivial canonical bundle
\jour Ann. of Math. \yr 1981 \vol 113 \pages 45--66
\endref
\ref \key Sch \by W. Schmid
\paper Variation of Hodge structure: the singularities of the period mapping
\jour Invent. Math. \yr 1973 \vol 22 \pages 211-319
\endref
\ref \key SzHs \by W. C. Hsiang and R. H. Szczarba
\paper On the tangent bundle of Grassman Manifold
\jour Amer. J. Math \yr 1964 \vol 86  \pages 698--704
\endref
\ref \key Ti \by Tian, G.
\paper Smoothness of the Universal Deformation Space of Compact
Calabi-Yau Manifolds and its Weil-Petersson Metric
\inbook Mathematical Aspects of String Theory
\eds S. T. Yau
\publ World Scientific \publaddr Singapore \pages 629--646
\endref
\ref \key Tod \by A. Todorov
\paper The Weil-Petersson Geometry of the Moduli Space of $SU(n)$
(Calabi-Yau) Manifolds I
\jour Commun. Math Phys \yr 1989 \vol 126 \pages 325--346
\endref
\enddocument